\documentclass[english]{article}
\usepackage[english]{babel}
\usepackage{latexsym,amssymb,amsmath}

\newcounter{theorem}

\renewcommand{\thetheorem}{\arabic{theorem}}

\newcommand{\theor}{\par\refstepcounter{theorem}{\bf Theorem \thetheorem .}\,\,}

\begin{document}

\centerline{\huge\textbf{Representation of Quaternionic}}
\vspace{.5cm}
\centerline{\huge\textbf{$\psi$--hyperholomorphic Functions}}

\vspace{.75cm}
\centerline{\Large Tetiana Kuzmenko and Vitalii Shpakivskyi}

\begin{abstract}
In the algebra of complex quaternions $\mathbb{H(C)}$ we consider for the first time left-- and right--$\psi$--hyperholomorphic functions. 
We justify the transition in left-- and right--$\psi$--hyperholomorphic functions to a simpler basis  i.e. to the Cartan basis.
Using Cartan's basis we find the solution of Cauchy--Fueter equation. By the same method we find a representation of left--$\psi$--hyperholomorphic
function in a special case.


\end{abstract}

\section{Introduction}\label{sec1}

Our main object of interest is the set which  is  usually  called  the  set  of  complex  quaternions  and which  is
 traditionally  denoted  as $\mathbb{H(C)}$. It turns  out  to be an  associative, non--commutative  complex  algebra
 generated  by the elements   $1$,  $ I $,  $J$, $K$ such that the following multiplication rules  hold:
\begin{equation*}\label{asterix1}
\begin{split}
I^2 =J^2  & = K^2 = IJK=-1 , \\
I J = - JI = K , \quad JK & = -KJ = I , \quad KI = - IK = J.
\end{split}
\end{equation*}
For $\mathbb{H(C)}$ another name, the algebra of biquaternions, is  used   also.

Consider in $\mathbb{H(C)}$ another set $\{e_1,e_2,e_3,e_4\}$, which is Cartan's
basis \cite{Cartan} such that

\begin{equation}\label{def-e}
 e_1=\frac{1}{2}(1+i I), \quad e_2=\frac{1}{2}(1-i I), \quad
e_3=\frac{1}{2}(iJ - K), \quad e_4=\frac{1}{2}(iJ +K) ,
\end{equation}
 where $i$ is
the complex imaginary unit.   It  is  direct  to  check that  we  got  a  new  basis.

The multiplication table  can
be represented as
\begin{equation}\label{tabl}
\begin{tabular}{c||c|c|c|c|}
$\cdot$ & $e_1$ & $e_2$ & $e_3$ & $e_4$\\
\hline \hline
$e_1$ & $e_1$ & $0$ & $e_3$ & $0$\\
\hline
$e_2$ & $0$ & $e_2$ & $0$ & $e_4$\\
\hline
$e_3$ & $0$ & $e_3$ & $0$ & $e_1$\\
\hline
$e_4$ & $e_4$ & $0$ & $e_2$ & $0$\\
\hline
\end{tabular}\,\,.
\end{equation}
The unit $1$ can be decomposed as $1=e_1+e_2$\,.

Note that the subalgebra with the basis $\{e_1,e_2\} $ is the
algebra of bicomplex numbers $\mathbb{BC}$ or Segre's algebra of
commutative quaternions (see, e.g., \cite{Alpay,Shapiro}).

The following relations holds:
\begin{equation}\label{reverse}
1=e_1+e_2\,,\quad I=-ie_1+ie_2\,,\quad J=-ie_3-ie_4\,,\quad
K=e_4-e_3\,.
\end{equation}
Of course,  formulas  (\ref{def-e}) and  (\ref{reverse}),  give the transition from one basis to the  other.

\section{Left-- and right--$\psi$--hyperholomorphic
functions} 

Let $\psi_1\,,\psi_2\,,\psi_3\,,\psi_4$ be fixed elements in
$\mathbb{H(C)}$ with the following representations in the Cartan's
basis:
$$\psi_1:=\sum\limits_{s=1}^4\alpha_se_s\,,\quad \alpha_s \in \mathbb{C},\qquad
\psi_2:=\sum\limits_{s=1}^4\beta_se_s\,,\quad \beta_s \in
\mathbb{C},$$
$$\psi_3:=\sum\limits_{s=1}^4\gamma_se_s\,,\quad \gamma_s \in \mathbb{C},\qquad
\psi_4:=\sum\limits_{s=1}^4\delta_se_s\,,\quad \delta_s \in
\mathbb{C}.$$ Consider a variable
$z=z_1e_1+z_2e_2+z_3e_3+z_4e_4\,,\,z_s\in\mathbb{C},\,s=1,2,3,4$ and
consider a function
$$f(z)=\sum\limits_{s=1}^4f_s(z_1,z_2,z_3,z_4)e_s\,,\quad f_s:\Omega\rightarrow \mathbb{H(C)},$$
where $\Omega$ is a domain in $\mathbb{C}^4$. Let components $f_s$,
$s=1,2,3,4$, are holomorphic functions of four complex variables
$z_1,z_2,z_3,z_4$ in $\Omega$.

Consider the operators

\begin{equation}\label{psi-left-op}
  ^{\psi}D:=\psi_1\frac{\partial }{\partial z_1}+\psi_2\frac{\partial }{\partial z_2}+\psi_3\frac{\partial }{\partial z_3}
  +\psi_4\frac{\partial }{\partial z_4}\,,
\end{equation}
\begin{equation}\label{psi-right-op}
   D^\psi:=\frac{\partial }{\partial z_1}\,\psi_1+\frac{\partial }{\partial z_2}\,\psi_2+\frac{\partial}{\partial z_3}\,\psi_3
   +\frac{\partial }{\partial z_4}\,\psi_4\,.
\end{equation}

A function
$f:\Omega\rightarrow\mathbb{H(C)}$, $\Omega\subset\mathbb{C}^4$, is
called \emph{left--$\psi$--hyperholomorphic (or
right--$\psi$--hyperholomorphic)} if components $f_s$ are holomorphic
functions of four complex variables $z_1,z_2,z_3,z_4$ in $\Omega$,
and $f$ satisfies the equation
\begin{equation}\label{psi-left}
^{\psi}D[f](z)=0.
\end{equation}
(or $ D^{\psi}[f](z)=0. $)

The class of $\psi$--hyperholomorphic functions in the real quaternions algebra is introduced for the first time by M.~V.~Shapiro and N.~L.~Vasilevski in 
the papers \cite{Shapiro-Vasilevski-1,Shapiro-Vasilevski-2}.
Since then, these functions have attracted the attention of many researchers. K.~Gürlebeck and his student H.~M.Nguyen pay a special attention to the
 applications of $\psi$--hyperholomorphic functions. See, for example, the papers \cite{Gurl-Nguen-13,Gurl-Nguen-15,Bock-Gurl-Nguen-15} and dissertation of H.~M.~Nguyen
 \cite{Nguen}. We note also that operators (\ref{psi-left-op}) and (\ref{psi-right-op}) are also called the weighted Dirac operators.
 Analysis and application of such operators are studied in papers  \cite{Vanegas,Vanegas-CV}.
 
 There are different generalizations of $\psi$--hyperholomorphic functions, which are being actively researched. 
Recently, generalizations to the case of fractional derivatives have become interesting. We will mark the works 
\cite{Cervantes-23,Cervantes-23-1}. 

Also began to consider operators of a more general form than (\ref{psi-left-op}). Namely, in paper \cite{Cervantes-23-2} investigated an operator of the form
$$
^{\psi}_\alpha D[f]:=\alpha f+\psi_1\frac{\partial f}{\partial z_1}+\psi_2\frac{\partial f}{\partial z_2}+\psi_3\frac{\partial f}{\partial z_3}
  +\psi_4\frac{\partial f}{\partial z_4}\,.
$$

In the paper \cite{Blaya-Math-Nachrichten}  it is develop the theory of co-called  $(\phi,\psi)$--hyperholomorphic functions. Following a
matrix approach, for such functions a generalized Borel--Pompeiu formula and the corresponding Plemelj-Sokhotski formulae are
established. Research from paper \cite{Blaya-Math-Nachrichten} was continued in the papers \cite{Blaya-Comp-Meth,Santiesteban,Santiesteban-Bull,Serrano-Blaya}.

At the same time, the problem of representation (or description in the explicit form) of hyperholomorphic and $\psi$--hyperholomorphic functions is open. 
This paper is devoted to solving this problem.

\subsection{Examples}

At first, we consider examples of left-- and right--$\psi$--hyperholomorphic functions.

\textbf{Example 1.}
 Consider a domain $\Omega \subset \mathbb{C}^2 \simeq \mathbb{BC}$ and consider a variable $\zeta=z_1e_1+z_2e_2$, and a
function $f:\Omega\rightarrow \mathbb{H(C)}$ of the form
$$f=\sum\limits_{s=1}^4f_s(z_1,z_2)e_s\,,\quad f_s:\Omega\rightarrow\mathbb{C}.$$

This  should  be  understood as  follows.  We identify $\mathbb{C}^2$ and $\mathbb {BC}$  after  which  the set $\Omega$
  in  $\mathbb {BC}$  becomes  a  subset  in $\mathbb{H(C)}$,  not  in $\mathbb{C}^2$;  next
   we  consider  some  objects  as  being  situated  in  $\mathbb{H(C)}$.  In particular,
   the set  $\Omega$ is situated  in  $\mathbb{H(C)}$.  When saying  that  the domain of $f$
    is  in  $\mathbb{H(C)}$ we  mean  already  the previous identifications.  Hence  we
      work with functions  with  both  domains  and  ranges  in  $\mathbb{H(C)}$.
       Thus  $\zeta$  is  in  a  domain  in $\mathbb{H(C)}$:  we imbed everything  in   $\mathbb{H(C)}$.

With  these  agreements we introduce the following  definitions.

A function $f:\Omega\rightarrow\mathbb{H(C)},\,\Omega\subset\mathbb{BC},$ is called \emph{right-$\mathbb{BC}$-hyperholomorphic}
if there exists an element $f'_r(\zeta)$ such that
\begin{equation}\label{r-B}
  \lim\limits_{\varepsilon\rightarrow 0}\frac{f(\zeta+\varepsilon
  h)-f(\zeta)}{\varepsilon}=
h\cdot f'_r(\zeta)\qquad \forall\,h\in\mathbb{BC}.
\end{equation}

A function $f:\Omega\rightarrow\mathbb{H(C)},\,\Omega\subset\mathbb{BC},$ is called \emph{left-$\mathbb{BC}$-hyperholomorphic}
if there exists an element $f'_l(\zeta)$ such that
\begin{equation}\label{l-B}
  \lim\limits_{\varepsilon\rightarrow 0}\frac{f(\zeta+\varepsilon
  h)-f(\zeta)}{\varepsilon}=
f'_l(\zeta)\cdot h\qquad \forall\,h\in\mathbb{BC}.
\end{equation}

Condition (\ref{r-B}) implies
\begin{equation}\label{A1}
  \frac{\partial f}{\partial z_1}=e_1f'_r(\zeta)\quad \text{for}\quad h=e_1
\end{equation}
and
\begin{equation}\label{A2}
  \frac{\partial f}{\partial z_2}=e_2f'_r(\zeta)\quad \text{for}\quad h=e_2.
\end{equation}

From (\ref{A1}) and (\ref{A2}) follows the analog of the
Cauchy--Riemann condition
\begin{equation}\label{C-R-r}
  e_2\frac{\partial f}{\partial z_1}=e_1\frac{\partial f}{\partial z_2}.
\end{equation}
Analogously, from (\ref{l-B}) follows
\begin{equation}\label{C-R-l}
  \frac{\partial f}{\partial z_1}\,e_2=\frac{\partial f}{\partial z_2}\,e_1.
\end{equation}

Thus, right-- and left--$\mathbb{BC}$--hyperholomorphic function
generalize holomorphic function theory in algebra $\mathbb{BC}$ (see,
e.g., \cite{Alpay,Shapiro}).

It is easy to see that the set of right- and
left--$\mathbb{BC}$--hyperholomorphic functions is a subset of
 left--$\psi$--hyperholomorphic and right--$\psi$--hyperholomorphic
function, respectively. Indeed,
 for $\zeta=z_1e_1+z_2e_2$ the equality (\ref{C-R-r}) has the form of the equality (\ref{psi-left})
 with $\psi_1=e_2$, $\psi_2=-e_1$, $\psi_3=\psi_4=0$. Analogously, left--$\mathbb{BC}$--hyperholomorphic functions is a
 subset of a set of right--$\psi$--hyperholomorphic functions.

Another example of mappings from the domain in $\mathbb{R}^3$ into
the algebra $\mathbb{H(C)}$, which are a particular case of left-- and right--$\psi$--hyperholomorphic
functions, is considered in
\cite{Kuzmenko-Shpakiv-19,Kuzmenko-Shpakiv-18}.

\textbf{Example 2.}
In (\ref{psi-left}) we set
$\psi_1=1,\,\psi_2=I,\,\psi_3=J,\,\psi_4=K$. In this case
$$\alpha_1=\alpha_2=1,\quad \alpha_3=\alpha_4=0,\quad \beta_1=-i,\quad \beta_2=i,\quad \beta_3=\beta_4=0,$$
$$\gamma_1=\gamma_2=0,\quad \gamma_3=-i,\quad \gamma_4=-i,\quad \delta_1=\delta_2=0,\quad \delta_3=-1,\quad\delta_4=1.$$
Then (\ref{psi-left}) takes the form
$$\frac{\partial f}{\partial z_1}+I\frac{\partial f}{\partial z_2}+J\frac{\partial f}{\partial z_3}+K\frac{\partial f}{\partial z_4}=0$$
that is well-known Cauchy--Fueter type equation (see, e.g.,
\cite{Shapiro-2014,Reyes-2010}).

\subsection{Main property of left-- and right--$\psi$--hyperholomorphic
functions}

\theor\label{theor-1} 
\emph{The property of function $f$ to be
left--$\psi$--hyperholomorphic (or right--$\psi$--hyperholomorphic) does not
 depend on a basis in which are given $f$ and $\psi$.}

\textbf{Proof.} Let us  prove the theorem for the case  left--$\psi$--hyperholomorphic functions.
 Let $\{e_1,e_2,e_3,e_4\}$ be the Cartan basis in
$\mathbb{H(C)}$ and $\{i_1,i_2,i_3,i_4\}$ another basis in
$\mathbb{H(C)}$. It means that
$$
  \begin{aligned}
&e_1=k_1i_1+k_2i_2+k_3i_3+k_4i_4,\\
&e_2=m_1i_1+m_2i_2+m_3i_3+m_4i_4,\\
&e_3=n_1i_1+n_2i_2+n_3i_3+n_4i_4,\\
&e_4=r_1i_1+r_2i_2+r_3i_3+r_4i_4,
\end{aligned}
$$
where $k_i,m_i,n_i,r_i,\,\,i=1,2,3,4,$ are complex numbers.

Consider the equation
\begin{equation}\label{6}
  D[f](t):=e_1\frac{\partial f}{\partial t_1} + e_2\frac{\partial f}{\partial t_2}+
  e_3\frac{\partial f}{\partial t_3}+e_4\frac{\partial f}{\partial t_4}=0,
\end{equation}
where $t:=t_1e_1+t_2e_1+t_3e_3+t_4e_4$,\,\,\,
$t_1,t_2,t_3,t_4\in\mathbb{C}$. Now we passing in $t$ to Cartan
basis. Then
$$t=i_1(t_1k_1+t_2m_1t_3+n_1+t_4r_1)+i_2(t_1k_2+t_2m_2+t_3n_2+t_4r_2)$$
$$+i_3(t_1k_3+t_2m_3+t_3n_3+t_4r_3)+i_4(t_1k_4+t_2m_4+t_3n_4+t_4r_4).$$

We set
\begin{equation}\label{7++}
 \begin{aligned}
&z_1:=t_1k_1+t_2m_1+t_3n_1+t_4r_1,\\
&z_2:=t_1k_2+t_2m_2+t_3n_2+t_4r_2,\\
&z_3:=t_1k_3+t_2m_3+t_3n_3+t_4r_3,\\
&z_4:=t_1k_4+t_2m_4+t_3n_4+t_4r_4.
\end{aligned}
\end{equation}

From equalities (\ref{7++}) we obtain
$$
 \begin{aligned}
&\frac{\partial f}{\partial t_1}=k_1\frac{\partial f}{\partial z_1}+k_2\frac{\partial f}{\partial z_2}
+k_3\frac{\partial f}{\partial z_3}+k_4\frac{\partial f}{\partial z_4}\,,\\
&\frac{\partial f}{\partial t_2}=m_1\frac{\partial f}{\partial z_1}+m_2\frac{\partial f}{\partial z_2}
+m_3\frac{\partial f}{\partial z_3}+m_4\frac{\partial f}{\partial z_4}\,,\\
&\frac{\partial f}{\partial t_3}=n_1\frac{\partial f}{\partial z_1}+n_2\frac{\partial f}{\partial z_2}
+n_3\frac{\partial f}{\partial z_3}+n_4\frac{\partial f}{\partial z_4}\,,\\
&\frac{\partial f}{\partial t_4}=r_1\frac{\partial f}{\partial
z_1}+r_2\frac{\partial f}{\partial z_2}+r_3\frac{\partial
f}{\partial z_3}+r_4\frac{\partial f}{\partial z_4}\,.
\end{aligned}
$$

Then equation (\ref{6}) is equivalent to the following equation
$$D[f](t)=(e_1k_1+e_2m_1+e_3n_1+e_4r_1)\frac{\partial f}{\partial z_1}+(e_1k_2+e_2m_2+e_3n_2+e_4r_2)\frac{\partial f}{\partial z_2}$$
$$+(e_1k_3+e_2m_3+e_3n_3+e_4r_3)\frac{\partial f}{\partial z_3}+(e_1k_4+e_2m_4+e_3n_4+e_4r_4)\frac{\partial f}{\partial z_4}$$
$$=\Big[i_1(k_1^2+m_1^2+n_1^2+r_1^2)+i_2(k_1k_2+m_1m_2+n_1n_2+r_1r_2)$$
$$+i_3(k_1k_3+m_1m_3+n_1n_3+r_1r_3)+i_4(k_1k_4+m_1m_4+n_1n_4+r_1r_4)\Big]\frac{\partial f}{\partial z_1}$$
$$+\Big[i_1(k_1k_2+m_1m_2+n_1n_2+r_1r_2)+i_2(k_2^2+m_2^2+n_2^2+r_2^2)$$
$$+i_3(k_2k_3+m_2m_3+n_2n_3+r_2r_3)+i_4(k_2k_4+m_2m_4+n_2n_4+r_2r_4)\Big]\frac{\partial f}{\partial z_2}$$
$$+\Big[i_1(k_1k_3+m_1m_3+n_1n_3+r_1r_3)+i_2(k_2k_3+m_2m_3+n_2n_3+r_2r_3)$$
$$+i_3(k_3^2+m_3^2+n_3^2+r_3^2)+i_4(k_3k_4+m_3m_4+n_3n_4+r_3r_4)\Big]\frac{\partial f}{\partial z_3}$$
$$+\Big[i_1(k_1k_4+m_1m_4+n_1n_4+r_1r_4)+i_2(k_2k_4+m_2m_4+n_2n_4+r_2r_4)$$
$$+i_3(k_3k_4+m_3m_4+n_3n_4+r_3r_4)+i_4(k_4^2+m_4^2+n_4^2+r_4^2)\Big]\frac{\partial f}{\partial z_4}$$
$$=:\psi_1\frac{\partial f}{\partial z_1}+\psi_2\frac{\partial f}{\partial z_2}+
\psi_3\frac{\partial f}{\partial z_3}+\psi_4\frac{\partial f}{\partial z_4}=0.$$

From this Theorem follows that in future investigation it is enough to consider constants
 $\psi$ and function $f$ in the simplest basis, i.e. in Cartan basis.

\section{Application to solving Cauchy--Fueter type equation}

Now, we will establish a connection between solutions of the equation
\begin{equation}\label{6}
  D[f](t):=\frac{\partial f}{\partial t_0} + I\frac{\partial f}{\partial t_1}+
  J\frac{\partial f}{\partial t_2}+K\frac{\partial f}{\partial t_3}=0,
\end{equation}
where $t:=t_0+t_1I+t_2J+t_3K$,\,\,\, $t_0,t_1,t_2,t_3\in\mathbb{C}$,
and the solutions of equations (\ref{psi-left}). For this purpose,
in $t$ we passing to Cartan basis. We have
$$t=t_0(e_1+e_2)+t_1(-ie_1+ie_2)+t_2(-ie_3-ie_4)+t_3(e_4-e_3)$$
$$=(t_0-it_1)e_1+(t_0+it_1)e_2+(-it_2-t_3)e_3+(-it_2+t_3)e_4.$$

We set
\begin{equation}\label{7}
  z_1:=t_0-it_1,\quad z_2:=t_0+it_1,\quad z_3:=-it_2-t_3,\quad z_4:=-it_2+t_3.
\end{equation}

From equalities (\ref{7}) we obtain
$$\frac{\partial f}{\partial t_0}=\frac{\partial f}{\partial z_1}+\frac{\partial f}{\partial z_2},\qquad
\frac{\partial f}{\partial t_1}=-i\frac{\partial f}{\partial z_1}+i\frac{\partial f}{\partial z_2},$$
$$\frac{\partial f}{\partial t_2}=-i\frac{\partial f}{\partial z_3}-i\frac{\partial f}{\partial z_4},\qquad
\frac{\partial f}{\partial t_3}=-\frac{\partial f}{\partial z_3}+\frac{\partial f}{\partial z_4}.$$

Then equation (\ref{6}) is equivalent to the following equation
$$D[f]=\frac{\partial f}{\partial z_1}+\frac{\partial f}{\partial z_2}-iI\frac{\partial f}{\partial z_1}+
iI\frac{\partial f}{\partial z_2}-iJ\frac{\partial f}{\partial z_3}-iJ\frac{\partial f}{\partial z_4}-
\frac{\partial f}{\partial z_3}+K\frac{\partial f}{\partial z_4}$$
$$=(1-iI)\frac{\partial f}{\partial z_1}+(1+iI)\frac{\partial f}{\partial z_2}+(-iJ-K)\frac{\partial f}{\partial z_3}
+(-iJ+K)\frac{\partial f}{\partial z_4}$$
$$=2\left(e_2\frac{\partial f}{\partial z_1}+e_1\frac{\partial f}{\partial z_2}-e_4\frac{\partial f}{\partial z_3}-
e_3\frac{\partial f}{\partial z_4}\right)=0.$$
Thus, we proved the following theorem

\theor\label{theor-2} \emph{A function $f$ of the variable
$t=t_0+t_1I+t_2J+t_3K$ satisfies equation (\ref{6})
 if and only if the function $f$ of the variable $z=z_1e_1+z_2e_2+z_3e_3+z_4e_4$ satisfies the equation
\begin{equation}\label{theor-8}
  e_2\frac{\partial f}{\partial z_1}+e_1\frac{\partial f}{\partial z_2}-e_4\frac{\partial f}{\partial z_3}-e_3\frac{\partial f}{\partial z_4}=0,
\end{equation}
where $z$ and $t$ are related by equalities (\ref{7}).}

Now, we solve equation (\ref{theor-8}).
$$e_2\frac{\partial f}{\partial z_1}=e_2\left(\frac{\partial f_1}{\partial z_1}e_1+\frac{\partial f_2}{\partial z_1}e_2
+\frac{\partial f_3}{\partial z_1}e_3+\frac{\partial f_4}{\partial
z_1}e_4\right)=\frac{\partial f_2}{\partial z_1}e_2+\frac{\partial
f_4}{\partial z_1}e_4\,,$$
$$e_1\frac{\partial f}{\partial z_2}=\frac{\partial f_1}{\partial z_2}e_1+\frac{\partial f_3}{\partial z_2}e_3\,,$$
$$e_4\frac{\partial f}{\partial z_3}=\frac{\partial f_1}{\partial z_3}e_4+\frac{\partial f_3}{\partial z_3}e_2\,,$$
$$e_3\frac{\partial f}{\partial z_4}=\frac{\partial f_2}{\partial z_4}e_3+\frac{\partial f_4}{\partial z_4}e_1\,.$$
Then equation (\ref{theor-8}) is equivalent to the system
$$\frac{\partial f_1}{\partial z_2}=\frac{\partial f_4}{\partial z_4}\,,\qquad \frac{\partial f_2}{\partial z_1}=\frac{\partial f_3}{\partial z_3}\,,$$
$$\frac{\partial f_3}{\partial z_2}=\frac{\partial f_2}{\partial z_4}\,,\qquad \frac{\partial f_4}{\partial z_1}=\frac{\partial f_1}{\partial z_3}\,.$$

We have pair of systems
\begin{equation}\label{syst-9}
\frac{\partial f_1}{\partial z_2}=\frac{\partial f_4}{\partial
z_4}\,,\qquad \frac{\partial f_1}{\partial z_3}=\frac{\partial
f_4}{\partial z_1}
\end{equation}
and
\begin{equation}\label{syst-10}
\frac{\partial f_2}{\partial z_1}=\frac{\partial f_3}{\partial
z_3}\,,\qquad \frac{\partial f_2}{\partial z_4}=\frac{\partial
f_3}{\partial z_2}\,.
\end{equation}

A solution of system (\ref{syst-9}), in a simple connected domain
$\Omega$, is an arbitrary holomorphic function
$$f_1=f_1(z_2,z_3)$$ and
$$f_4=z_4\frac{\partial f_1}{\partial z_2}+z_1\frac{\partial f_1}{\partial z_3}.$$

A solution of system (\ref{syst-10}), in a simple connected domain
$\Omega$, is an arbitrary holomorphic function
$$f_2=f_2(z_1,z_4)$$
and
$$f_3=z_3\frac{\partial f_2}{\partial z_1}+z_2\frac{\partial f_2}{\partial z_4}.$$

Thus, we have the following solution of equation (\ref{theor-8}):
$$ f(z)=f_1(z_2,z_3)e_1+f_2(z_1,z_4)e_2$$
\begin{equation}\label{11}
 +\left(z_3\frac{\partial f_2}{\partial z_1}+z_2\frac{\partial f_2}{\partial z_4}\right)e_3
 +\left(z_4\frac{\partial f_1}{\partial z_2}+z_1\frac{\partial f_1}{\partial z_3}\right)e_4\,.
\end{equation}

Thus, accordingly to Theorem \ref{theor-2} we obtain

\theor
 \emph{In a simple connected domain, function (\ref{11})
 in which $z_1,z_2,z_3,z_4$ are given by relations (\ref{7}),
 satisfies  equation (\ref{6}).}

\textbf{Proposition 1.}
  In a simple connected domain,  function (\ref{11})
 satisfies the four-dimensional complex Laplace equation
\begin{equation}\label{cor-2}
  \Delta_{\mathbb{C}^4}f:=\frac{\partial^2 f}{\partial z_1^2}+\frac{\partial^2 f}{\partial z_2^2}
+\frac{\partial^2 f}{\partial z_3^2}+\frac{\partial^2 f}{\partial z_4^2}=0.
\end{equation}

About equation (\ref{cor-2}) and its relation with the Cauchy-Fueter
equation see in \cite{Shapiro-2014}.

\section{Representation of left--$\psi$--hyperholomorphic
function in a special case}

Now we will find a general solution of equation (\ref{psi-left}) for
a special choice of parameters $\psi_1,\psi_2,\psi_3$ and $\psi_4$.
For this purpose, we reduce equation (\ref{psi-left}) to a system of
four PDEs. We have
$$\psi_1\frac{\partial f}{\partial z_1}=\left(\alpha_1e_1+\alpha_2e_2+\alpha_3e_3+\alpha_4e_4\right)
\left(\frac{\partial f_1}{\partial z_1}e_1+\frac{\partial f_2}{\partial z_1}e_2+\frac{\partial f_3}{\partial z_1}e_3+
\frac{\partial f_4}{\partial z_1}e_4\right)$$
$$=\frac{\partial f_1}{\partial z_1}\alpha_1e_1+\frac{\partial f_3}{\partial z_1}\alpha_1e_3+
\frac{\partial f_2}{\partial z_1}\alpha_2e_2+\frac{\partial f_4}{\partial z_1}\alpha_2e_4$$
$$+\frac{\partial f_2}{\partial z_1}\alpha_3e_3+\frac{\partial f_4}{\partial z_1}\alpha_3e_1+
\frac{\partial f_1}{\partial z_1}\alpha_4e_4+\frac{\partial f_3}{\partial z_1}\alpha_4e_2$$
$$=\frac{\partial}{\partial z_1}\left(\alpha_1f_1+\alpha_3f_4\right)e_1+\frac{\partial}{\partial z_1}\left(\alpha_2f_2+\alpha_4f_3\right)e_2$$
$$+\frac{\partial}{\partial z_1}\left(\alpha_1f_3+\alpha_3f_2\right)e_3+\frac{\partial}{\partial z_1}\left(\alpha_2f_4+\alpha_4f_1\right)e_4.$$

Similarly
$$\psi_2\frac{\partial f}{\partial z_2}=\frac{\partial}{\partial z_2}\left(\beta_1f_1+\beta_3f_4\right)e_1
+\frac{\partial}{\partial z_2}\left(\beta_2f_2+\beta_4f_3\right)e_2$$
$$+\frac{\partial}{\partial z_2}\left(\beta_1f_3+\beta_3f_2\right)e_3+\frac{\partial}{\partial z_2}\left(\beta_2f_4+\beta_4f_1\right)e_4,$$

$$\psi_3\frac{\partial f}{\partial z_3}=\frac{\partial}{\partial z_3}\left(\gamma_1f_1+\gamma_3f_4\right)e_1+
\frac{\partial}{\partial z_3}\left(\gamma_2f_2+\gamma_4f_3\right)e_2$$
$$+\frac{\partial}{\partial z_3}\left(\gamma_1f_3+\gamma_3f_2\right)e_3+\frac{\partial}{\partial z_3}\left(\gamma_2f_4+\gamma_4f_1\right)e_4,$$

$$\psi_4\frac{\partial f}{\partial z_4}=\frac{\partial}{\partial z_4}\left(\delta_1f_1+\delta_3f_4\right)e_1
+\frac{\partial}{\partial z_4}\left(\delta_2f_2+\delta_4f_3\right)e_2$$
$$+\frac{\partial}{\partial z_4}\left(\delta_1f_3+\delta_3f_2\right)e_3+\frac{\partial}{\partial z_4}\left(\delta_2f_4+\delta_4f_1\right)e_4.$$

Then equation (\ref{psi-left}) is equivalent to the following system
$$
\frac{\partial}{\partial
z_1}(\alpha_1f_1+\alpha_3f_4)+\frac{\partial}{\partial
z_2}(\beta_1f_1+\beta_3f_4) +\frac{\partial}{\partial
z_3}(\gamma_1f_1+\gamma_3f_4)+\frac{\partial}{\partial
z_4}(\delta_1f_1+\delta_3f_4)=0,$$ $$ \frac{\partial}{\partial
z_1}(\alpha_2f_2+\alpha_4f_3)+\frac{\partial}{\partial
z_2}(\beta_2f_2+\beta_4f_3) +\frac{\partial}{\partial
z_3}(\gamma_2f_2+\gamma_4f_3)+\frac{\partial}{\partial
z_4}(\delta_2f_2+\delta_4f_3)=0,$$
\begin{equation}\label{syst-12}\vspace{-2mm}
\end{equation}
$$\frac{\partial}{\partial
z_1}(\alpha_1f_3+\alpha_3f_2)+\frac{\partial}{\partial
z_2}(\beta_1f_3+\beta_3f_2) +\frac{\partial}{\partial
z_3}(\gamma_1f_3+\gamma_3f_2)+\frac{\partial}{\partial
z_4}(\delta_1f_3+\delta_3f_2)=0,$$
$$\frac{\partial}{\partial
z_1}(\alpha_2f_4+\alpha_4f_1)+\frac{\partial}{\partial
z_2}(\beta_2f_4+\beta_4f_1) +\frac{\partial}{\partial
z_3}(\gamma_2f_4+\gamma_4f_1)+\frac{\partial}{\partial
z_4}(\delta_2f_4+\delta_4f_1)=0.$$ \vspace{2mm}

\theor
\label{theor-3} \emph{For
\begin{equation}\label{13}
\begin{aligned}
&\psi_1=\alpha_1e_1+\alpha_2e_2+\alpha_3e_3+\alpha_4e_4,\qquad \alpha_1\alpha_2\neq\alpha_3\alpha_4\,,\\
&\psi_2=\lambda\alpha_1e_1+\mu\alpha_2e_2+\mu\alpha_3e_3+\lambda\alpha_4e_4,\\
&\psi_3=\theta\alpha_1e_1+\vartheta\alpha_2e_2+\vartheta\alpha_3e_3+\theta\alpha_4e_4,\\
&\psi_4=\nu\alpha_1e_1+\eta\alpha_2e_2+\eta\alpha_3e_3+\nu\alpha_4e_4,
\end{aligned}
\end{equation}
where $\alpha_1,\alpha_2,\alpha_3,\alpha_4,\lambda,\mu,\theta,\vartheta,\nu,\eta$ are an arbitrary complex
numbers, the general solution of equation (\ref{psi-left}) is of the form
\begin{equation}\label{13-G}
f(z)=f_1(\widetilde{\zeta}_2,\widetilde{\zeta}_3,\widetilde{\zeta}_4)e_1+f_2(\zeta_2,\zeta_3,\zeta_4)e_2+
f_3(\widetilde{\zeta}_2,\widetilde{\zeta}_3,\widetilde{\zeta}_4)e_3+f_4(\zeta_2,\zeta_3,\zeta_4)e_4,
\end{equation}
where
\begin{equation}\label{13-1}
\begin{aligned}
  &\widetilde{\zeta}_2:=\lambda z_1-z_2,\quad \widetilde{\zeta}_3:=\theta z_1-z_3,\quad \widetilde{\zeta}_4:=\nu z_1-z_4,\\
  &\zeta_2:=\mu z_1-z_2,\quad \zeta_3:=\vartheta z_1-z_3,\quad \zeta_4:=\eta z_1-z_4.
\end{aligned}
\end{equation}}

\textbf{Proof.} For given parameters (\ref{13}) the first equation
of  system (\ref{syst-12}) takes the form
$$  \frac{\partial}{\partial z_1}\left(\alpha_1f_1+\alpha_3f_4\right)+\frac{\partial}{\partial z_2}\left(\lambda\alpha_1f_1+\mu\alpha_3f_4\right)+$$
\begin{equation}\label{14}
+\frac{\partial}{\partial z_3}\left(\theta\alpha_1f_1+\vartheta\alpha_3f_4\right)+\frac{\partial}{\partial z_4}\left(\nu\alpha_1f_1+\eta\alpha_3f_4\right)=0.
\end{equation}

Similarly, for given parameters (\ref{13}) the fourth equation of
system (\ref{syst-12}) takes the form
$$\frac{\partial}{\partial z_1}\left(\alpha_4f_1+\alpha_2f_4\right)+\frac{\partial}{\partial z_2}\left(\lambda\alpha_4f_1+\mu\alpha_2f_4\right)$$
\begin{equation}\label{15}
  +\frac{\partial}{\partial z_3}\left(\theta\alpha_4f_1+\vartheta\alpha_2f_4\right)+\frac{\partial}{\partial z_4}\left(\nu\alpha_4f_1+\eta\alpha_2f_4\right)=0.
\end{equation}

Consider the difference between equation (\ref{14}) multiplied by $\alpha_2$ and equation (\ref{15}) multiplied by $\alpha_3$.
Then we obtain the following equation
$$\frac{\partial}{\partial z_1}\Big(f_1(\alpha_1\alpha_2-\alpha_3\alpha_4)+f_4(\alpha_2\alpha_3-\alpha_2\alpha_3)\Big)$$
$$+\frac{\partial}{\partial z_2}\Big(f_1(\lambda\alpha_1\alpha_2-\lambda\alpha_3\alpha_4)+f_4(\mu\alpha_2\alpha_3-\mu\alpha_2\alpha_3)\Big)$$
$$+\frac{\partial}{\partial z_3}\Big(f_1(\theta\alpha_1\alpha_2-\theta\alpha_3\alpha_4)+f_4(\vartheta\alpha_2\alpha_3-\vartheta\alpha_2\alpha_3)\Big)$$
$$+\frac{\partial}{\partial z_4}\Big(f_1(\nu\alpha_1\alpha_2-\nu\alpha_3\alpha_4)+f_4(\eta\alpha_2\alpha_3-\eta\alpha_2\alpha_3)\Big)=0.$$

Thus, we obtain the equation
\begin{equation}\label{16}
  \frac{\partial f_1}{\partial z_1}+\lambda\frac{\partial f_1}{\partial z_2}+\theta\frac{\partial f_1}{\partial z_3}+\nu\frac{\partial f_1}{\partial z_4}=0.
\end{equation}

For equation (\ref{16}) consider the characteristic equation
\begin{equation}\label{17}
  \frac{dz_1}{1}=\frac{dz_2}{\lambda}=\frac{dz_3}{\theta}=\frac{dz_4}{\nu}.
\end{equation}

The solutions of system (\ref{17}) are the following integrals
$$c_2=\lambda z_1-z_2,\quad c_3=\theta z_1-z_3,\quad c_4=\nu z_1-z_4.$$

Therefore, the general solution of equation (\ref{16}) has the form
$$f_1=f_1(\widetilde{\zeta}_2,\widetilde{\zeta}_3,\widetilde{\zeta}_4),$$
where $\widetilde{\zeta}_2,\widetilde{\zeta}_3,\widetilde{\zeta}_4$ are defined by equalities (\ref{13-1}).

Note that polynomials (\ref{13-1}) are similarly to the well-known Fueter's polynomials \cite{Alpay-2005}.

Similarly we obtain the representations for the components $f_2,f_3,f_4$.

Thus, formula (\ref{13-G}) given representation of every
left--$\psi$--hyperholomorphic function.

\textbf{Tetiana Kuzmenko}\\
Department of Fundamental Sciences\\
Zhytomyr Military Institute\\
Prospect Myru, 22\\
10004, Zhytomyr, Ukraine\\
kuzmenko.ts15@gmail.com

\vspace{1cm}
\textbf{Vitalii Shpakivskyi}\\
Department of Complex Analysis and Potential Theory\\
Institute of Mathematics of the National Academy of Sciences of Ukraine\\
3, Tereschenkivska st.\\
01024, Kyiv-4, Ukraine\\
shpakivskyi86@gmail.com

\end{document}